\documentclass[12pt,a4paper]{article}
\usepackage{amssymb,amsmath,amsthm}
\usepackage{graphicx}
\usepackage{tikz}
\usetikzlibrary{patterns}
\usetikzlibrary{arrows}
\usetikzlibrary{decorations.markings}
\usetikzlibrary{arrows}
\usepackage{hyperref}
\usepackage{xcolor}

\newcommand\ex{\ensuremath{\mathrm{ex}}}
\newcommand\exo{\ensuremath{\mathrm{ex_o}}}
\newcommand\exd{\ensuremath{\mathrm{ex_d}}}

\theoremstyle{plain}
\newtheorem{theorem}{Theorem}[section]
\newtheorem{lemma}[theorem]{Lemma}

\newtheorem{proposition}[theorem]{Proposition}

\theoremstyle{definition}

\newtheorem{definition}[theorem]{Definition}

\newtheorem{claim}[theorem]{Claim}

\newtheorem*{lemma*}{Lemma}
\newtheorem*{thm*}{Theorem}

\newcommand\cref[1]{Corollary~\ref{cor:#1}}

\title{Stars and paths in oriented graphs}

\author{Xuanrui Hu$^{1}$, Hanyi Wang$^{2}$, Yuefang Sun$^{3}$, Bo Deng$^{4}$
\\
$^{1}$ School of Mathematics and Statistics, Ningbo University,\\
Ningbo 315211, China, huxuanrui1108@163.com\\
$^{2}$ School of Mathematics and Statistics, Ningbo University,\\
Ningbo 315211, China, wanghanyi119@163.com\\
$^{3}$ Corresponding author. School of Mathematics and Statistics,\\
Ningbo University,
Ningbo 315211, China, sunyuefang@nbu.edu.cn\\
$^{4}$ School of Mathematics and Statistics, Qinghai Normal University,\\
Xining 810016, China,
State Key Laboratory of Tibetan Intelligence,\\
Xining, Qinghai, 810008, China, dengbo450@163.com}
\date{}

\begin{document}

\maketitle
\begin{center}
Dedicated to Professor Gyula O.H. Katona on the occasion of his 85th birthday.
\end{center}
\begin{abstract}
Let $\overrightarrow{F}$ be a fixed oriented graph. Its oriented Tur\'{a}n number, written $\exo(n,\overrightarrow{F})$, is the largest size of an $n$-vertex oriented graph with no copy of $\overrightarrow{F}$. This parameter extends the classical notion of Tur\'{a}n number to oriented graphs.
In this paper, we investigate the supersaturation results for the directed path $P_k~(k\geq 2)$ and all oriented stars $\overrightarrow{S}_{p,q}~(0\leq p\leq q$),  and all oriented paths of order at most four.
\end{abstract}
\vspace{0.3cm}


\vspace{0.3cm}

\section{Introduction}
All digraphs considered in this paper are simple; in particular, we exclude loops and parallel arcs from consideration. The notation $|S|$ is adopted for the cardinality of a set $S$. We use $K_k$ to denote the complete graph of order $k$, and  $T(n,k)$ represents the $\textit{Tur\'{a}n graph}$, which is an $n$-vertex complete $k$-partite graph with partition sets balanced so that any two parts differ in size by at most one.
 
Given a digraph $D$, its vertex and arc sets are represented by $V(D)$ and $A(D)$. For a vertex $v$, its \textit{total degree} $d_D(v)$ is defined as the sum of out-degree $d^+_D(v)$ and in-degree $d^-_D(v)$ in $D$. Let $N_D^+(v)$ and $N_D^-(v)$ denote the out-neighbourhood and in-neighbourhood of $v$. A digraph with no directed cycles is called \textit{acyclic}. A digraph containing no $2$-cycles is called an \textit{oriented graph}. An oriented graph is \textit{antidirected} if all its vertices are either sources or sinks; equivalently, it is bipartite with all arcs oriented from one partition class to the other. We abbreviate a tournament of order $n$ as an $n$-tournament. We call an $n$-tournament $T$ \textit{transitive} if the following holds for all distinct vertices $u_1,u_2,u_3$: Whenever $T$ contains arcs $\overrightarrow{u_1u_2}$ and $\overrightarrow{u_2u_3}$, it contains an arc from $u_1$ to $u_3$. It can be checked that every transitive $n$-tournament possesses a vertex enumeration $u_1,u_2,\ldots, u_n$ such that every arc goes from $u_i$ to $u_j$ for $1\leq i<j\leq n$. Consequently, this is the unique $n$-vertex acyclic tournament, denoted $\overrightarrow{TT}_n$. We write $\overrightarrow{AT}_3$ for the antidirected path on four vertices and $\overrightarrow{OT}_3$ for the oriented path $u_1\rightarrow u_2\rightarrow u_3\leftarrow u_4$. The notation $N(\overrightarrow{H},\overrightarrow{G})$ gives the total number of $\overrightarrow{H}$-copies embedded in $\overrightarrow{G}$. If $N(\overrightarrow{H},\overrightarrow{G})=0$, then $\overrightarrow{G}$ is $\overrightarrow{H}$-\textit{free}.

\begin{definition}[Oriented Tur\'an number]Given an oriented graph $\overrightarrow{F}$, the \textit{oriented Tur\'an number} $\exo(n,\overrightarrow{F})$ is the largest size among all $n$-vertex oriented graphs that are $\overrightarrow{F}$-free. 
\end{definition}
As an oriented counterpart of the classical Tur\'an number $\ex(n,F)$, the oriented Tur\'an number $\exo(n,\overrightarrow{F})$ describes the extremal arc number for $\overrightarrow{F}$-free oriented graphs. By contrast, $\exd(n,H)$ is a related parameter, defined as the largest size of an $H$-free digraph of order $n$.
Any oriented graph satisfies $$\exo(n,\overrightarrow{F})\le \exd(n,\overrightarrow{F})\le 2\exo(n,\overrightarrow{F}).$$ We note that $\exd(n,H)$ has been studied extensively, starting from the work of Brown and Harary \cite{Brown-Harary}.

The well-known Erdős-Stone-Simonovits Theorem~\cite{Erdos-Stone, Erdos-Simonovits} asserts that  $$\ex(n,F)=\left(\frac{\chi(F)-2}{\chi(F)-1}+o(1)\right)\binom{n}{2}$$ holds for any graph $F$. It completely characterizes the asymptotics of $\ex(n,F)$ provided $\chi(F)\geq 3$. 
The study of $\exo(n,\overrightarrow{F})$ was first undertaken by Valadkhan~\cite{vala}, who derived an Erdős-Stone-Simonovits-type Theorem~\cite{Erdos-Stone, Erdos-Simonovits} for oriented graphs. In this oriented version, the compressibility $z(\overrightarrow{F})$ takes the place of the chromatic number $\chi(F)$ appearing in the original theorem.

P. Erd\H{o}s and M. Simonovits~\cite{Erdos-Simonovits} stated the general Erdős–Simonovits supersaturation theorem.
\begin{theorem}[P. Erd\H{o}s and M. Simonovits~\cite{Erdos-Simonovits}]\label{PM}
    Given $\varepsilon>0$ and an $h$-vertex graph $H$, there exist a positive constant \(\delta=\delta(H,\varepsilon)\) and a positive integer \(n_0\) such that every \(n\)-vertex graph \(G\) with \(n\ge n_0\) and  $|E(G)|>\ex(n,H)+\varepsilon n^2$ satisfies
    $$N(H,G)\geq \delta n^h.$$
\end{theorem}
X. Hu and Y. Sun~\cite{Hu-Sun} established the existence of the oriented Tur\'{a}n density $\pi_o(\overrightarrow{F})$. This result is an oriented analogue of the classical Katona-Nemetz-Simonovits Theorem~\cite{Katona-Nemetz-Simonovits}. Based on this, they established the oriented version of Theorem~\ref{PM} and investigated the supersaturation of $\overrightarrow{TT}_3$, showing that every $n$-vertex oriented graph $\overrightarrow{G}$ with $|A(\overrightarrow{G})|=|E(T(n,3))|+1$ satisfies $N(\overrightarrow{TT}_3,\overrightarrow{G})\geq \left\lfloor\frac{2n}{3}\right\rfloor.$

 \begin{theorem}[X. Hu and Y. Sun~\cite{Hu-Sun}]\label{superori}
Given $\varepsilon >0$ and an $h$-vertex oriented graph $\overrightarrow{F}$, there exist $\delta>0$ and a positive integer $n_0$ such that every $n$-vertex oriented graph $\overrightarrow{G}$ with $n\geq n_0$ and $|A(\overrightarrow{G})|\geq
\left(\pi_o(\overrightarrow{F})+\varepsilon\right)\binom{n}{2}$ satisfies
$$N(\overrightarrow{F},\overrightarrow{G})\geq \delta \binom{n}{h}.$$
\end{theorem}
 Several upper bounds on the extremal number of antidirected graphs have been established. Since $\exo(n,\overrightarrow{F})\leq \exd(n,\overrightarrow{F})$, any upper bound for $\exd(n,\overrightarrow{F})$ also yields the corresponding bound for $\exo(n,\overrightarrow{F})$. For an antidirected tree $\overrightarrow{F}$ with $k\geq 2$ arcs, Graham \cite{gra} proved that $\exd(n,\overrightarrow{F})=O(n)$. Burr \cite{burr} subsequently obtained the explicit estimate $\exd(n,\overrightarrow{F})\leq 4kn$. For antidirected paths, Klimo\u{s}ov\'{a} and Stein \cite{ks} established that $\exo(n,\overrightarrow{F})\leq (3k-4)n/2$ whenever $k\geq 3$. More recently, Grzesik and Skrzypczyk \cite{gs} improved this estimate by showing that
\[
\exo(n,\overrightarrow{F})
\leq (k-1+\sqrt{k-3})n
\]
for every antidirected path of length $k\geq 4$. 
Gerbner, Hu and Sun~\cite{Gerbner-Hu-Sun} extended the F\"{u}redi--Alon--Krivelevich--Sudakov Theorem~\cite{Furedi, Alon-Krivelevich-Sudakov}, which generalizes the K{\H{o}}v{\'a}ri--S\'os--Tur\'an theorem \cite{kst} to antidirected graphs, obtaining the following result.

\begin{theorem}[Gerbner-Hu-Sun~\cite{Gerbner-Hu-Sun}]\label{orientedFAKS}
Suppose $\overrightarrow{H}$ is obtained from a bipartite graph with parts $A$ and $B$ by directing each edge from $A$ to $B$, and $d_{\overrightarrow{H}}^+(u)\le r$ for all $u\in A$. Then for some constant $c=c(\overrightarrow{H})>0$, we have
$$\exo(n,\overrightarrow{H})\le cn^{2-\frac{1}{r}}.$$
 \end{theorem}

Gerbner, Hu and Sun~\cite{Gerbner-Hu-Sun} obtained exact values of the oriented Tur\'{a}n numbers for several classes of oriented graphs, including $\overrightarrow{P}_k$, three distinct orientations of the path on four vertices, and oriented stars.
We recall the following results: 
\begin{proposition}\label{Pk}
    $\exo(n,\overrightarrow{P}_k)=|E(T(n,k-1))|.$
\end{proposition}

\begin{proposition}\label{3-path1}
     For sufficiently large $n$, $\exo(n,\overrightarrow{OT}_3)=\left\lfloor \dfrac{n^2}{4} \right\rfloor$.
\end{proposition}

\begin{proposition}\label{3-path2}
   For the oriented graph $\overrightarrow{F}$ on four vertices $v_1,v_2,v_3,v_4$ with arc set $\{\overrightarrow{v_2v_1},\overrightarrow{v_2v_3},\overrightarrow{v_3v_4}\}$, we have
$$\exo(n,\overrightarrow{F})=\left\lfloor \dfrac{n^2}{4} \right\rfloor$$
for sufficiently large $n$.
\end{proposition}


\begin{proposition}\label{antidirected}
     $\exo(n,\overrightarrow{AT}_3)=2n-3$ $(n>1)$.
\end{proposition}

\begin{proposition}\label{Spq}
$\exo(n,\overrightarrow{S}_{p,q})=\begin{cases}
    (p-1)n+\left\lfloor\frac{(n+q-p)^{2}}{4}\right\rfloor,& 1\le p\le q~(\text{for~large}~n)\\
    ({q}-1){n},&{p}=0\\
    \end{cases}$.
\end{proposition}
Obviously, the oriented 3-paths in Proposition~\ref{3-path2} are equivalent to those in Proposition~\ref{3-path1} by symmetry: This follows by reversing all arcs and relabeling the vertices.

R\'edei\cite{red} proved that every $k$-tournament contains a $\overrightarrow{P}_k$.
\begin{theorem}\label{path k}
    Every $k$-tournament contains a $\overrightarrow{P}_k$.
\end{theorem}
Using Theorem~\ref{path k} and Proposition~\ref{Pk}, we can derive
supersaturation bounds for $\overrightarrow{P}_k$ from
corresponding bounds for $K_k$ in the underlying graph.

\subsection{Our results}
In this paper, we investigate the supersaturation results for the directed path $P_k~(k\geq 2)$ and all oriented stars $\overrightarrow{S}_{p,q}~(0\leq p\leq q$),  and all oriented paths of order at most four. The core supersaturation phenomenon we consider is the following: if an oriented graph $\overrightarrow{G}$ has number of oriented edges above $\exo(n,\overrightarrow{F})$, then it has ``many" copies of $\overrightarrow{F}$. 

Gerbner, Hu and Sun~\cite{Gerbner-Hu-Sun} proved that 
$\exo(n,\overrightarrow{P}_k)=|E(T(n,k-1))|$. 
We establish the following supersaturation theorem for $\overrightarrow{P}_k$.

\begin{theorem}\label{P4}
For sufficiently large $n$, write $n=(k-1)s+t$,
where $s\ge1$ and $0\le t\le k-2$ are integers.
Every oriented graph $\overrightarrow{G}$ on $n$ vertices
with exactly $|E(T(n,k-1))|+1$ arcs satisfies
\[
N(\overrightarrow{P}_k,\overrightarrow{G})
\ge
\begin{cases}
s^{k-2}, & t=0,\\[2mm]
s^{k-1-t}(s+1)^{t-1}, & 1\le t\le k-2.
\end{cases}
\]
\end{theorem}
Next, by Proposition~\ref{3-path1}, we obtain two supersaturation results for \(\overrightarrow{OT}_3\): an $n$-vertex oriented graph with \(\left\lfloor \frac{n^2}{4}\right\rfloor+1\) arcs contains at least \(\Omega(n^2)\) copies of \(\overrightarrow{OT}_3\), and one with \(tn^2\) arcs (\(t>\tfrac14\)) contains at least \(\Omega(n^4)\) copies of \(\overrightarrow{OT}_3\). By symmetry from Proposition \ref{3-path2} , we can get the same supersaturation results on $\overrightarrow{F}$.
\begin{theorem}\label{ot3}
   Let $\overrightarrow{G}$ be an $n$-vertex oriented graph, where $n$ is sufficiently large. Then the following assertions hold:
    \begin{description}
  \item[(i)]  If $|A(\overrightarrow{G})|= \left\lfloor\frac{n^2}{4}\right\rfloor+1$, then
    $$N(\overrightarrow{OT}_3, \overrightarrow{G})\geq \left\lceil \frac{n}{2}\right\rceil\left(\left\lfloor\frac{n}{2}\right\rfloor-2\right).$$
    \item[(ii)] If $|A(\overrightarrow{G})| = t n^2$ for some $t > \frac{1}{4}$, then
\[
N(\overrightarrow{OT}_3,\overrightarrow{G})
\ge
\frac{t^2(4t-1)^2}{9}n^4
-\frac{1}{2}n^3.
\]
    \end{description}
\end{theorem}

We then turn to $\overrightarrow{AT}_3$. By Proposition~\ref{antidirected},  $\exo(n,\overrightarrow{AT}_3)=2n-3$ ($n>1$). We obtain the following supersaturation result of $\overrightarrow{AT}_3$:

\begin{theorem}\label{antidirected 4-path}
     Let $\overrightarrow{G}$ be an $n$-vertex oriented graph, where $n$ is sufficiently large. If $|A(\overrightarrow{G})|= tn$ for some  $t\geq 3$, then
    $$N(\overrightarrow{AT}_3, \overrightarrow{G})\geq tn\left(2t-3-\sqrt{2t-3}\right).$$
\end{theorem}

Finally, using Proposition \ref{Spq}, we obtain two supersaturation results for oriented stars.

\begin{theorem}\label{spq}
    Fix positive integers $p$ and $q$ with $ p\leq q$, and let $\overrightarrow{G}$ be an $n$-vertex oriented graph, where $n$ is sufficiently large. If $|A(\overrightarrow{G})|=(p-1)n+\left\lfloor \frac{(n+q-p)^2}{4}\right\rfloor +1$, then
$$N(\overrightarrow{S}_{p,q},\overrightarrow{G})\geq \binom{\lfloor\frac{n+p-q}{2}\rfloor}{p}.$$
\end{theorem}
\begin{theorem}\label{S0p}
    Fix a positive integers $q$, and let $\overrightarrow{G}$ be an $n$-vertex oriented graph, where $n$ is sufficiently large. If $|A(\overrightarrow{G})|=(q-1)n+k$ for some integer $k\geq 1$, then
$$N(\overrightarrow{S}_{0,q},\overrightarrow{G})\geq k.$$
\end{theorem}

\section{Proof of Theorem~\ref{P4}}

\vskip 0.5cm
\noindent
\begin{lemma}\label{critical}\cite{red}
    Let $G$ be a graph with sufficiently large order $n$. If $|E(G)|=|E(T(n,k-1))|+1$, then
    $$N(K_k,G)\geq \prod_{i=0}^{k-3}\left\lfloor{\frac{n+i}{k-1}}\right\rfloor.$$
    
\end{lemma}
\noindent
{\bf Theorem~\ref{P4}:} {\em   For sufficiently large $n$, write $n=(k-1)s+t$,
where $s\ge1$ and $0\le t\le k-2$ are integers.
Every oriented graph $\overrightarrow{G}$ on $n$ vertices
with exactly $|E(T(n,k-1))|+1$ arcs satisfies
\[
N(\overrightarrow{P}_k,\overrightarrow{G})
\ge
\begin{cases}
s^{k-2}, & t=0,\\[2mm]
s^{k-1-t}(s+1)^{t-1}, & 1\le t\le k-2.
\end{cases}
\]
    }

\begin{proof}
Let $G$ denote the underlying graph of $\overrightarrow{G}$. Since
\[
|E(G)|=|A(\overrightarrow{G})| =|E(T(n,k-1))|+1,
\]
Lemma~\ref{critical} implies that
$$N(K_k,G)
\geq \prod_{i=0}^{k-3}\left\lfloor{\frac{n+i}{k-1}}\right\rfloor.$$
Each copy of $K_k$ in $G$ induces a $k$-tournament in $\overrightarrow{G}$. By Theorem~\ref{path k}, every $k$-tournament contains a copy of $\overrightarrow{P}_k$.\\
Therefore, 
$$N(\overrightarrow{P}_k,\overrightarrow{G})=N(K_k,G) \geq \prod_{i=0}^{k-3}\left\lfloor{\frac{n+i}{k-1}}\right\rfloor =
\begin{cases}
s^{k-2},& t=0\\
       s^{k-1-t}(s+1)^{t-1},&1\leq t\leq k-2\\
\end{cases}.$$
This finish the proof.
\end{proof}

\section{Proof of Theorem~\ref{ot3}}
{\bf Theorem~\ref{ot3}:} {\em   Let $\overrightarrow{G}$ be an $n$-vertex oriented graph, where $n$ is sufficiently large. Then the following assertions hold:
    \begin{description}
  \item[(i)]  If $|A(\overrightarrow{G})|= \left\lfloor\frac{n^2}{4}\right\rfloor+1$, then
    $$N(\overrightarrow{OT}_3, \overrightarrow{G})\geq \left\lceil \frac{n}{2}\right\rceil\left(\left\lfloor\frac{n}{2}\right\rfloor-2\right).$$
    \item[(ii)] If $|A(\overrightarrow{G})| = t n^2$ for some $t > \frac{1}{4}$, then
\[
N(\overrightarrow{OT}_3,\overrightarrow{G})
\ge
\frac{t^2(4t-1)^2}{9}n^4
-\frac{1}{2}n^3.
\]
    \end{description}
    }

 \begin{proof}
 \noindent{\bf Part (i)}
    If \(N( \overrightarrow{OT}_3,\overrightarrow G)>n^2/4\), the desired bound holds trivially. We thus assume $N(\overrightarrow{OT}_3, \overrightarrow{G})\leq \frac{n^2}{4}$. Define
  $$S_1=\{u: d^-_{\overrightarrow{G}}(u)<\frac{n}{10}\},~~~S_2=V(\overrightarrow{G})\setminus S_1.$$
  Since
  $$|A(\overrightarrow{G})|=\sum_{u\in V(\overrightarrow{G})}d^-_{\overrightarrow{G}}(u)=\left\lfloor\frac{n^2}{4}\right\rfloor+1\leq |S_1|\frac{n}{10}+|S_2|n=\frac{n^2}{10}+\frac{9n|S_2|}{10},$$
We obtain
\begin{equation*}
     \begin{aligned}\label{b1}
     |S_2|\geq \frac{n}{6}.
\end{aligned}
 \tag{3.1}
\end{equation*}
For an arc $\overrightarrow{uv}\in A(\overrightarrow{G})$ with $v\in S_2$, we use
$N(\overrightarrow{OT}_3,\overrightarrow{uv})$
to count those $\overrightarrow{OT}_3$-copies in which
$\overrightarrow{uv}$ serves as the middle arc. Thus, 
$$N(\overrightarrow{OT}_3,\overrightarrow{uv})\geq d^-_{\overrightarrow{G}}(u)(d^-_{\overrightarrow{G}}(v)-2)\geq d^-_{\overrightarrow{G}}(u)(\frac{n}{10}-2).$$
Summing over all arcs terminating in $S_2$, 
$$\frac{n^2}{4}\geq N(\overrightarrow{OT}_3, \overrightarrow{G})\geq (\frac{n}{10}-2)\sum_{u\in V(\overrightarrow{G})}d^-_{\overrightarrow{G}}(u)d^+_{S_2}(u)$$
Consequently, 
\begin{equation*}
     \begin{aligned}\label{b2}
\sum_{u\in V(\overrightarrow{G})}d^-_{\overrightarrow{G}}(u)d^+_{S_2}(u)\leq \frac{5n^2}{2n-40}.
\end{aligned}
 \tag{3.2}
\end{equation*}
For disjoint vertex sets $X,Y$, let $e(X,Y)$ denote the number of arcs directed from $X$ to $Y$ and let $e(\overrightarrow{G}[X])$ denote the number of arcs in the induced subgraph $\overrightarrow{G}[X]$. An arc is called a \textit{bad arc} if it is not directed from $S_1$ to $S_2$. For $u\in S_1$ and $v\in S_2$, we call $\overrightarrow{uv}$ a \textit{missing arc} if $\overrightarrow{uv}\notin A(\overrightarrow{G})$. Set $$p=\sum_{u\in S_1}d^-_{\overrightarrow{G}}(u)+e(\overrightarrow{G}[S_2]),~~~q=|S_1||S_2|-e(S_1,S_2).$$
Thus, $p$ and $q$ are the numbers of bad arcs and missing arcs, respectively.
Since 
$$|A(\overrightarrow{G})|=\left\lfloor\frac{n^2}{4}\right\rfloor+1=e(S_1,S_2)+\sum_{u\in S_1}d^-_{\overrightarrow{G}}(u)+e(\overrightarrow{G}[S_2])=|S_1||S_2|-q+p,$$
we have 
\begin{equation*}
     \begin{aligned}\label{b3}
     p-q=\left\lfloor\frac{n^2}{4}\right\rfloor+1-|S_1||S_2|\geq 1,
\end{aligned}
 \tag{3.3}
\end{equation*}
where the inequality follows from $|S_1||S_2|\leq \left\lfloor\frac{n^2}{4}\right\rfloor.$ \\
We now expand the sum:
\begin{equation*}
     \begin{aligned}\label{b4}
\sum_{u\in V(\overrightarrow{G})}d^-_{\overrightarrow{G}}(u)d^+_{S_2}(u)&=\sum_{u\in S_1}d^-_{\overrightarrow{G}}(u)d^+_{S_2}(u)+\sum_{u\in S_2}d^-_{\overrightarrow{G}}(u)d^+_{S_2}(u)\\
&\geq \sum_{u\in S_1}d^-_{\overrightarrow{G}}(u)d^+_{S_2}(u)+\frac{n}{10}e({\overrightarrow{G}}[S_2])\\
&=|S_2|\sum_{u\in S_1}d^-_{\overrightarrow{G}}(u)-\sum_{u\in S_1}d^-_{\overrightarrow{G}}(u)\left(|S_2|-d^+_{S_2}(u)\right)+\frac{n}{10}e({\overrightarrow{G}}[S_2])\\
&\geq |S_2|\sum_{u\in S_1}d^-_{\overrightarrow{G}}(u)-\frac{n}{10}q+\frac{n}{10}e({\overrightarrow{G}}[S_2])\\
&=|S_2|\sum_{u\in S_1}d^-_{\overrightarrow{G}}(u)-\frac{n}{10}\left(q-p+\sum_{u\in S_1}d^-_{\overrightarrow{G}}(u)\right)\\
&=\left(|S_2|-\frac{n}{10}\right)\sum_{u\in S_1}d^-_{\overrightarrow{G}}(u)-\frac{n}{10}(q-p)\\
&\geq \left(|S_2|-\frac{n}{10}\right)\sum_{u\in S_1}d^-_{\overrightarrow{G}}(u),
\end{aligned}
\tag{3.4}
\end{equation*}
where the second inequality follows from the definition of $S_1$ and $\sum_{u\in S_1}\left(|S_2|-d_{S_2}^+(u)\right)=q.$\\
By \eqref{b1}, $|S_2|-n/ 10\geq n/ 15$.
Together with \eqref{b2}, this yields
$$\sum_{u\in S_1}d^-_{\overrightarrow{G}}(u)=O(1).$$
Since also 
$$\sum_{u\in V(\overrightarrow{G})}d^-_{\overrightarrow{G}}(u)d^+_{S_2}(u)\geq \sum_{u\in S_2}d^-_{\overrightarrow{G}}(u)d^+_{S_2}(u)\geq \frac{n}{10}e({\overrightarrow{G}}[S_2]),$$
by \eqref{b2}, we get $e({\overrightarrow{G}}[S_2])=O(1).$
Hence, 
\begin{equation*}
     \begin{aligned}
     p=\sum_{u\in S_1}d^-_{\overrightarrow{G}}(u)+e({\overrightarrow{G}}[S_2])=O(1).
\end{aligned}
\end{equation*}
It now follows from \eqref{b3} that
\begin{equation*}
     \begin{aligned}\label{b6}
    q\leq p-1=O(1), ~~~\left\lfloor\frac{n^2}{4}\right\rfloor-|S_1||S_2|=p-q-1=O(1). 
\end{aligned}
 \tag{3.5}
\end{equation*}
Note that $$(|S_1|-|S_2|)^2=n^2-4|S_1||S_2|=n^2-4\left\lfloor\frac{n^2}{4}\right\rfloor+4\left(\left\lfloor\frac{n^2}{4}\right\rfloor-|S_1||S_2|\right)=O(1),$$
which implies $$|S_1|=\frac{n}{2}+O(1),~~~|S_2|=\frac{n}{2}+O(1).$$
  Fix a bad arc $\overrightarrow{u_1v_1}\in A(\overrightarrow{G})$, and count copies of $\overrightarrow{OT}_3$ containing $\overrightarrow{u_1v_1}$ whose remaining two arcs are directed from $S_1$ to $S_2$. There are three cases.\\
CASE 1: $u_1,v_1\in S_1$. Any such copy takes the form \\
$$u_1\rightarrow v_1\rightarrow u_2\leftarrow v_2,~~~u_2\in S_2,~v_2\in S_1\setminus\{u_1,v_1\}.$$
By \eqref{b6}, there are \(q=O(1)\) missing arcs from \(S_1\) to \(S_2\). Each missing arc excludes at most $n$
candidate copies, so at most $qn=O(n)$ candidates are
excluded. Therefore, the number of such copies of $\overrightarrow{OT}_3$ is at least $$ |S_2|(|S_1|-2)-O(n)=\frac{n^2}{4}-O(n).$$
CASE 2: $u_1,v_1\in S_2$. Any such copy takes the form \\
$$u_2\rightarrow u_1\rightarrow v_1\leftarrow v_2,~~~u_2,~v_2\in S_1,~u_2\neq v_2.$$
Similarly, the number of such copies of $\overrightarrow{OT}_3$ is at least $$ |S_1|(|S_1|-1)-O(n)=\frac{n^2}{4}-O(n).$$
CASE 3: $u_1\in S_2,~v_1\in S_1$. Any such copy takes the form \\
$$u_1\rightarrow v_1\rightarrow u_2\leftarrow v_2,~~~u_2\in S_2\setminus\{u_1\},~v_2\in S_1\setminus\{v_1\}.$$
Similarly, the number of such copies of $\overrightarrow{OT}_3$ is at least $$ (|S_2|-1)(|S_1|-1)-O(n)=\frac{n^2}{4}-O(n).$$
Each counted copy contains exactly one bad arc, so the families from distinct bad arcs are pairwise disjoint. By definition of $p$, there are $p=O(1)$ bad arcs. If $p\geq 2$, then
$$N(\overrightarrow{OT}_3,\overrightarrow{G})\geq 2\left(\frac{n^2}{4}-O(n)\right)>\frac{n^2}{4},$$
a contradiction. Hence, $p\leq 1$. From \eqref{b6}, $p=1,~q=0,~|S_1||S_2|=\left\lfloor \frac{n^2}{4}\right\rfloor$.\\
Therefore, if $N(\overrightarrow{OT}_3,\overrightarrow{G})\leq \frac{n^2}{4}$, then $S_1$ and $S_2$ form a balanced partition, and all arcs directed from $S_1$ to $S_2$ are present. The unique bad arc $e$ is contained in either $\overrightarrow{G}[S_1]$ or $\overrightarrow{G}[S_2]$.\\
Thus,
\[
N(\overrightarrow{OT}_3,\overrightarrow{G})
=
\begin{cases}
|S_2|(|S_1|-2), & e\in A(\overrightarrow{G}[S_1])\\
|S_1|(|S_1|-1), & e\in A(\overrightarrow{G}[S_2])
\end{cases}.
\]
Since $\{|S_1|,|S_2|\}=\{\lfloor n/2\rfloor,\lceil n/2\rceil\}$,
we obtain
$$N(\overrightarrow{OT}_3,\overrightarrow{G})\geq \left\lceil \frac{n}{2}\right\rceil\left(\left\lfloor \frac{n}{2}\right\rfloor-2\right).$$

\noindent{\bf Part (ii)} 
Let $\overrightarrow G$ be an $n$-vertex oriented graph with $tn^2$ arcs for some $t>\frac{1}{4}$, where $n$ is sufficiently large.
We first establish the following claim.
  \begin{claim}\label{2-path} 
   \[
  N( \overrightarrow{P}_3,\overrightarrow G) \ge \dfrac{t(4t-1)}{3}n^3.
\]
\end{claim}
\begin{proof}[Proof of Claim \ref{2-path}]
 Let \(G\) be the underlying graph of $\overrightarrow{G}$.\\
Every $C_3$ in \(G\) corresponds to either a \(\overrightarrow{TT}_3\) or a \(\overrightarrow{C}_3\) in \(\overrightarrow{G}\), containing one and three copies of \(\overrightarrow{P}_3\), respectively.
Since copies arising from distinct triangles have distinct vertex sets, we obtain
\[
 N( \overrightarrow{P}_3,\overrightarrow G) \ge N(C_3,G).
\]
For any two adjacent vertices $u,v\in G$, we obtain
\[
\begin{aligned}
|N_G(u)\cap N_G(v)|
&=d_G(u)+d_G(v)-|N_G(u)\cup N_G(v)|\\
&\ge d_G(u)+d_G(v)-n.
\end{aligned}
\]
Since each $C_3$ is counted once for each of its three edges,
\[
\begin{aligned}
N(C_3,G)
&=\frac{1}{3}\left(\sum_{uv\in E(G)} |N_G(u)\cap N_G(v)|\right)\\
&\ge \frac{1}{3}\left(\sum_{uv\in E(G)}\bigl(d_G(u)+d_G(v)-n\bigr)\right)\\
&=\frac{1}{3}\left(\sum_{v\in V(G)}d_G(v)^2-tn^3\right)\\
&\geq \frac{1}{3n}\left(\sum_{v\in V(G)} d_G(v)\right)^2-\frac{t}{3}n^3\\
&=\frac{t(4t-1)}{3}n^3,
\end{aligned}
\]
where the last inequality follows from the Cauchy--Schwarz inequality, and the last equality uses $\sum_{v\in V(G)}d_G(v)=2tn^2$. This finish the proof of Claim~\ref{2-path}.
\end{proof}
We use $N(\overrightarrow{P}_3,v)$
to count those $\overrightarrow{P}_3$-copies terminating at $v$. 
Then
\[
N(\overrightarrow{P}_3,\overrightarrow{G})=\sum_{v\in V(\overrightarrow{G})}N(\overrightarrow{P}_3,v).
\]
We count copies of $\overrightarrow{OT}_3$ by extending directed paths $w_1\to u\to v$. For each such path, choose $w_2\in N_{\overrightarrow{G}}^-(v)$. This gives a copy of $\overrightarrow{OT}_3$ unless $w_2=u$ or $w_2=w_1$, and each copy arises uniquely in this way. The choices with $w_2=u$ are counted by $N(\overrightarrow{P}_3,\overrightarrow{G})$, whereas those with
$w_2=w_1$ are in one-to-one correspondence with copies of
$\overrightarrow{TT}_3$. Consequently,
\[
N(\overrightarrow{OT}_3,\overrightarrow{G})
=
\sum_{v\in V(\overrightarrow{G})}d^-_{\overrightarrow{G}}(v)N(\overrightarrow{P}_3,v)-\left(N(\overrightarrow P_3,\overrightarrow G)
  +N(\overrightarrow{TT}_3,\overrightarrow G)\right).
\]
Each vertex triple contributes at most $3$ to $N(\overrightarrow P_3,\overrightarrow G)
  +N(\overrightarrow{TT}_3,\overrightarrow G)$: a $\overrightarrow{C}_3$ contributes $3$, a $\overrightarrow{TT}_3$ contributes $2$, and all other triples contribute at most $1$. Thus,
\[
N(\overrightarrow{P}_3,\overrightarrow{G})+N(\overrightarrow{TT}_3,\overrightarrow{G})
\le 3\binom{n}{3}
\le \frac{n^3}{2}.
\]
Moreover,
$$N(\overrightarrow{P}_3,v)=\sum_{u\in N_{\overrightarrow{G}}^-(v)}d^-_{\overrightarrow{G}}(u)\le \sum_{u\in N^-_{\overrightarrow{G}}(v)} n
=nd^-_{\overrightarrow{G}}(v).$$
Therefore, by the Cauchy-Schwarz inequality and Claim~\ref{2-path}
\begin{align*}
N(\overrightarrow{OT}_3,\overrightarrow{G})
&=\sum_{v\in V(\overrightarrow{G})}d^-_{\overrightarrow{G}}(v)N(\overrightarrow{P}_3,v)-\left(N(\overrightarrow P_3,\overrightarrow G)
  +N(\overrightarrow{TT}_3,\overrightarrow G)\right)\\
&\ge \frac{1}{n}\sum_{v\in V(\overrightarrow{G})}N^2(\overrightarrow{P}_3,v)-\frac{1}{2}n^3\\
&\ge \frac{1}{n^2}\left(\sum_{v\in V(\overrightarrow{G})}N(\overrightarrow{P}_3,v)\right)^2-\frac{1}{2}n^3\\
&=\frac{N^2(\overrightarrow{P}_3,\overrightarrow{G})}{n^2}-\frac{1}{2}n^3\\
&\geq \frac{t^2(4t-1)^2}{9}n^4 - \frac{1}{2}n^3.
\end{align*}
This finsh the proof.

\end{proof}

\section{Proof of Theorem~\ref{antidirected 4-path}}
\begin{definition}
    For any three distinct vertices $u_1$, $u_2$ and $u_3$, given a $\overrightarrow{TT}_3$ constructed as
    $$u_3\leftarrow u_1\rightarrow u_2\leftarrow u_3,$$
    we call $\overrightarrow{u_1u_2}$ the \textit{middle arc} of this $\overrightarrow{TT}_3$. 
\end{definition}

\begin{definition}
    For any four distinct vertices $u$, $v$, $w_1$ and $w_2$, given a $\overrightarrow{AT}_3$ constructed as 
    $$w_1\leftarrow u\rightarrow v\leftarrow w_2,$$
    we call $\overrightarrow{uv}$ the \textit{middle arc} of this $\overrightarrow{AT}_3$.
\end{definition}
In the following, let $\overrightarrow{G}$ be an $n$-vertex oriented graph, where $n$ is sufficiently large, and 
let $C(\overrightarrow{G})=N(\overrightarrow{AT}_3,\overrightarrow{G})+N(\overrightarrow{TT}_3,\overrightarrow{G})$. 
For any $\overrightarrow{uv}\in A(\overrightarrow{G})$, let $N(\overrightarrow{AT}_3, \overrightarrow{uv})$ denote the number of copies of $\overrightarrow{AT}_3$ whose middle arc is $\overrightarrow{uv}$, let $N(\overrightarrow{TT}_3, \overrightarrow{uv})$ denote the number of copies of $\overrightarrow{TT}_3$ whose middle arc is $\overrightarrow{uv}$ and $C(\overrightarrow{uv})=N(\overrightarrow{AT}_3, \overrightarrow{uv})+N(\overrightarrow{TT}_3, \overrightarrow{uv})$. For every arc $\overrightarrow{uv}\in A(\overrightarrow{G})$, each choice of $w_1\in N^+(u)\backslash \{v\}$ and $w_2\in N^-(v)\backslash \{u\}$ (where $w_1$ and $w_2$ may coincide) yields exactly one copy of either $\overrightarrow{AT}_3$ or $\overrightarrow{TT}_3$. Therefore, \begin{equation*}
    \begin{aligned}
C(\overrightarrow{G})&=N(\overrightarrow{AT}_3,\overrightarrow{G})+N(\overrightarrow{TT}_3,\overrightarrow{G})\\&=\sum_{\overrightarrow{uv}\in A(\overrightarrow{G})}\left(d^+_{\overrightarrow{G}}(u)-1\right)\left(d^-_{\overrightarrow{G}}(v)-1\right)\\&=\sum_{u\in V(\overrightarrow{G})}\left(d^+_{\overrightarrow{G}}(u)-1\right)\left[\sum_{v\in N^+(u)}(d^-_{\overrightarrow{G}}(v)-1)\right].
\end{aligned}
\end{equation*}
{\bf Theorem~\ref{antidirected 4-path}:} {\em    Let $\overrightarrow{G}$ be an $n$-vertex oriented graph, where $n$ is sufficiently large. If $|A(\overrightarrow{G})|= tn$ for some  $t\geq 3$, then
    $$N(\overrightarrow{AT}_3, \overrightarrow{G})\geq tn\left(2t-3-\sqrt{2t-3}\right).$$}
\begin{proof}
Let $|A(\overrightarrow{G})|=tn(t\geq 3)$. Define two vertex sets 
$$U=\{u\in V(\overrightarrow{G}):d^+_{\overrightarrow{G}}(u)=1\},~W=\{v\in V(\overrightarrow{G}):d^-_{\overrightarrow{G}}(v)=1\},$$  
let $A(U,W)$ denote the set of arcs from $U$ to $W$. Obviously, an arc $\overrightarrow{uv}$ contributes no copy of either $\overrightarrow{AT}_3$ or $\overrightarrow{TT}_3$ as a middle arc if and only if $d_{\overrightarrow{G}}^+(u)=1$ or $d^-_{\overrightarrow{G}}(v)=1$.  Consequently, the total number of such arcs is $q=|U|+|W|-|A(U,W)|$. Note that $|U|,|W|<n$; otherwise, $|A(\overrightarrow{G})|=n$, contradicting
$|A(\overrightarrow{G})|=tn$ with $t\ge3$.\\
Let $$A_1=\{\overrightarrow{uv}\in A(\overrightarrow{G}): d^+_{\overrightarrow{G}}(u)\geq2,~d^-_{\overrightarrow{G}}(v)\geq2\},$$ then $|A_1|=|A(\overrightarrow{G})|-q.$
Let $\overrightarrow{G}_1$ be the arc-induced subgraph of $\overrightarrow{G}$ induced by $A_1$. Summing over all arcs in $A_1$, we obtain $$\sum_{u\in V(\overrightarrow{G}_1)}d^+_{\overrightarrow{G}_1}(u)=\sum_{v\in V(\overrightarrow{G}_1)}d^-_{\overrightarrow{G}_1}(v)=|A_1|.$$ 
By the definition of  $\overrightarrow{G}_1$, we have
$$d^+_{\overrightarrow{G}}(u)\ge d^+_{\overrightarrow{G}_1}(u),~ d^-_{\overrightarrow{G}}(v)\ge d^-_{\overrightarrow{G}_1}(v) . $$
In $\overrightarrow{G}$, for any arc $\overrightarrow{uv}\in A_1,$
 \begin{equation*}
    \begin{aligned}
C(\overrightarrow{uv})&=\left(d^+_{\overrightarrow{G}}(u)-1\right)\left(d^-_{\overrightarrow{G}}(v)-1\right)\\
&=\left(d^+_{\overrightarrow{G}}(u)-2\right)\left(d^-_{\overrightarrow{G}}(v)-2\right)+d^+_{\overrightarrow{G}}(u)+d^-_{\overrightarrow{G}}(v)-3\\
&\geq d^+_{\overrightarrow{G}}(u)+d^-_{\overrightarrow{G}}(v)-3\\
&\geq d^+_{\overrightarrow{G}_1}(u)+d^-_{\overrightarrow{G}_1}(v)-3.
\end{aligned}
\end{equation*}
 Summing over all arcs in $A_1$, we obtain \begin{equation}\label{eq:myeq1}
C(\overrightarrow{G})\geq \sum_{\overrightarrow{uv}\in A_1}\left[d^+_{\overrightarrow{G}_1}(u)+d^-_{\overrightarrow{G}_1}(v)-3\right]
=\sum_{u\in V(\overrightarrow{G}_1)}\left(d^+_{\overrightarrow{G}_1}(u)\right)^2
+\sum_{v\in V(\overrightarrow{G}_1)}\left(d^-_{\overrightarrow{G}_1}(v)\right)^2-3|A_1|.
\tag{4.1}
\end{equation}
 Note that at most $n-|U|$ vertices have positive out-degrees in $\overrightarrow{G}_1$ and at most $n-|W|$ vertices have positive in-degrees in $\overrightarrow{G}_1$. By the Cauchy–Schwarz inequality, $$\sum_{u\in V(\overrightarrow{G}_1)}\left({d^+_{\overrightarrow{G}_1}(u)}\right)^2
 \geq\frac{\left(\sum_{u\in V(\overrightarrow{G}_1)}d^+_{\overrightarrow{G}_1}(u)\right)^2}{\sum_{u\in V(\overrightarrow{G}_1)}1^2}\geq
 \frac{|A_1|^2}{n-|U|},$$ $$~\sum_{v\in V(\overrightarrow{G}_1)}\left({d^-_{\overrightarrow{G}_1}(v)}\right)^2
 \geq\frac{\left(\sum_{v\in V(\overrightarrow{G}_1)}d^-_{\overrightarrow{G}_1}(v)\right)^2}{\sum_{v\in V(\overrightarrow{G}_1)}1^2}\geq
 \frac{|A_1|^2}{n-|W|}.$$   
 Therefore, by \eqref{eq:myeq1}, 
 \begin{equation}\label{eq:myeq2}
C(\overrightarrow{G})\geq |A_1|^2\left(\frac{1}{n-|U|}+\frac{1}{n-|W|}\right)-3|A_1|.
\tag{4.2}
\end{equation} 
 Since $n-|U|>0$ and $n-|W|>0,$ 
 \begin{equation}\label{eq:myeq3}
\frac{1}{n-|U|}+\frac{1}{n-|W|}=\frac{n-|U|+n-|W|}{(n-|U|)(n-|W|)}\geq \frac{4}{n-|U|+n-|W|}=\frac{4}{2n-|U|-|W|},
\tag{4.3}
\end{equation}
 the inequality holds by $$\frac{n-|U|+n-|W|}{(n-|U|)(n-|W|)}-\frac{4}{n-|U|+n-|W|}=\frac{(n-|U|-n+|W|)^2}{(n-|U|+n-|W|)(n-|U|)(n-|W|)}\geq0.$$ 
 Since $|U|+|W|=q+|A(U,W)|\geq q$, we have $2n-|U|-|W|\leq 2n-q$. Therefore, by \eqref{eq:myeq2}, \eqref{eq:myeq3},
 \begin{equation}\label{eq:4}
    \begin{aligned}
C(\overrightarrow{G})&\geq |A_1|^2\left(\frac{1}{n-|U|}+\frac{1}{n-|W|}\right)-3|A_1|\\
&\geq \frac{4|A_1|^2}{2n-|U|-|W|}-3|A_1|\\
&\geq\frac{4|A_1|^2}{2n-q}-3|A_1|\\
&=\frac{4(tn-q)^2}{2n-q}-3(tn-q)\\
&\geq t(2t-3)n,
\end{aligned}\tag{4.4}
\end{equation}
where the last inequality holds since 
 $$\frac{4(tn-q)^2}{2n-q}-3(tn-q)-t(2t-3)n=\frac{q(2n(t-1)(t-3)+q)}{2n-q}\geq 0$$
 for $t\geq 3.$\\
 For each arc $\overrightarrow{uv}\in A(\overrightarrow{G})$,  $|N^+(u)\cap N^-(v)|\leq\min\{d^+(u)-1,d^-(v)-1\}$. Thus,  $$N(\overrightarrow{TT}_3,\overrightarrow{uv})\leq \min\{d^+(u)-1,d^-(v)-1\}\leq \sqrt{(d^+(u)-1)(d^-(v)-1)}=\sqrt{C(\overrightarrow{uv})}.$$
 Summing over all arcs and applying the Cauchy–Schwarz inequality,
 \begin{equation}\label{eq:5}
    \begin{aligned}
N(\overrightarrow{TT}_3,\overrightarrow{G})&\leq \sum_{\overrightarrow{uv}\in A(\overrightarrow{G})}\sqrt{C(\overrightarrow{uv})}\\
&=\sqrt{\left(\sum_{\overrightarrow{uv}\in{A(\overrightarrow{G})}}1\cdot \sqrt{C(\overrightarrow{uv})} \right)^2}\\
 &\leq \sqrt{\left(\sum_{\overrightarrow{uv}\in A(\overrightarrow{G})}1^2\right)\cdot \left(\sum_{\overrightarrow{uv}\in A(\overrightarrow{G})}\left(\sqrt{C(\overrightarrow{uv})}\right)^2\right)}\\
&=\sqrt{|A(\overrightarrow{G})|\sum_{\overrightarrow{uv}\in A(\overrightarrow{G})}C(\overrightarrow{uv})}\\
&=\sqrt{tnC(\overrightarrow{G})}. 
\end{aligned}\tag{4.5}
\end{equation}
Therefore, by \eqref{eq:4}, \eqref{eq:5}
 \begin{equation*}
    \begin{aligned}
N(\overrightarrow{AT}_3,\overrightarrow{G})&= C(\overrightarrow{G})-N(\overrightarrow{TT}_3,\overrightarrow{G})\\
&\geq C(\overrightarrow{G})-\sqrt{tnC(\overrightarrow{G})}\\
&\geq t(2t-3)n-\sqrt{tn\cdot t(2t-3)n}\\
&=tn\left(2t-3-\sqrt{2t-3}\right),
\end{aligned}
\end{equation*}
where the last inequality follows from the fact that the function $f(x)=x-\sqrt{tnx}$ is increasing for $x\geq t(2t-3)n$.\\
This finish the proof.   
\end{proof}

\section{Proof of Theorem~\ref{spq}}

{\bf Theorem~\ref{spq}:} {\em     Fix positive integers $p$ and $q$ with $ p\leq q$, and let $\overrightarrow{G}$ be an $n$-vertex oriented graph, where $n$ is sufficiently large. If $|A(\overrightarrow{G})|=(p-1)n+\left\lfloor \frac{(n+q-p)^2}{4}\right\rfloor +1$, then
$$N(\overrightarrow{S}_{p,q},\overrightarrow{G})\geq \binom{\lfloor\frac{n+p-q}{2}\rfloor}{p}.$$}

\begin{proof}
     Define $F(n)=(p-1)n+\left\lfloor \frac{(n+q-p)^2}{4}\right\rfloor=\exo(n,\overrightarrow{S}_{p,q}).$\\
    By definition of the oriented Tur\'{a}n number, every $\overrightarrow{S}_{p,q}$-free oriented graph $\overrightarrow{H}$ on $m$ vertices satisfies
 $$|A(\overrightarrow{H})|\leq F(m)=(p-1)m+\left\lfloor\frac{(m+q-p)^2}{4}\right\rfloor.$$
 Let $Z=\{v\in V(\overrightarrow{G}):~d_{\overrightarrow{G}}^-(v)\geq p,~d_{\overrightarrow{G}}^+(v)\geq q\}$. Consider the induced oriented graph $\overrightarrow{H}=\overrightarrow{G}-Z$.\\ 
 Clearly, $\overrightarrow{H}$ contains no $\overrightarrow{S}_{p,q}$, so 
 $|A(\overrightarrow{H})|\leq F(n-|Z|)$.\\
 Every arc not in $\overrightarrow{H}$ is incident to at least one vertex of $Z$, hence
  \begin{equation}\label{a2}
 \sum_{v\in Z}d_{\overrightarrow{G}}(v)\geq |A(\overrightarrow{G})|-|A(\overrightarrow{H})|\geq F(n)+1-F(n-|Z|).
 \tag{5.1}
\end{equation}
 For each $v\in Z$, we use $N(\overrightarrow{S}_{p,q},v)$ to count those $\overrightarrow{S}_{p,q}$-copies in which vertex $v$ serves as the center. Hence,
 $$N(\overrightarrow{S}_{p,q},v)=\binom{d^-_{\overrightarrow{G}}(v)}{p}\binom{d^+_{\overrightarrow{G}}(v)}{q}.$$
 We first establish the following claim.
 \begin{claim}\label{Pascal}
     For integers $a\geq p$, $b\geq q$ with $p\leq q$, we have
     $$\binom{a}{p}\binom{b}{q}\geq\binom{a+b-q}{p}.$$
 \end{claim}
\begin{proof}[Proof of Claim \ref{Pascal}]
    Set $r=b-q\geq 0.$ It suffices to show $$\binom{a}{p}\binom{q+r}{q}\geq \binom{a+r}{p}.$$
We observe that
\begin{equation*}
     \begin{aligned}
 \frac{\binom{a+r}{p}}{\binom{a}{p}}&=\prod_{i=1}^r\frac{a+i}{a-p+i}
 =\prod_{i=1}^r(1+\frac{p}{a-p+i})\\
 &\leq \prod_{i=1}^r(1+\frac{p}{i})
 \leq \prod_{i=1}^r(1+\frac{q}{i})\\
 &=\frac{(q+1)(q+2)\ldots (q+r)}{r!}=\binom{q+r}{q},
  \end{aligned}
\end{equation*}
Rearranging yields the desired inequality.
\end{proof}

By the definition of $Z$, we obtain
\begin{equation*}
     \begin{aligned}\label{a4}
N(\overrightarrow{S}_{p,q},\overrightarrow{G})&=\sum_{v\in Z}N(\overrightarrow{S}_{p,q},v)\\
&=\sum_{v\in Z}\binom{d_{\overrightarrow{G}}^-(v)}{p}\binom{d_{\overrightarrow{G}}^+(v)}{q}\\
&\geq \sum_{v\in Z}\binom{d_{\overrightarrow{G}}(v)-q}{p}\\
&\geq |Z|\binom{\lfloor\frac{1}{|Z|}\sum_{v\in Z}(d_{\overrightarrow{G}}(v)-q)\rfloor}{p},
 \end{aligned}
 \tag{5.2}
\end{equation*}
where the last inequality follows from Jensen's inequality. We divide the proof into three cases according to $|Z|$.\\
CASE 1: $|Z|=1$. Let $Z=\{v\}$, by \eqref{a2},
\begin{equation*}
     \begin{aligned}
d_{\overrightarrow{G}}(v)-q&\geq F(n)+1-F(n-1)-q\\
&=(p-1)n+\left\lfloor \frac{(n+q-p)^2}{4}\right\rfloor+1-(p-1)(n-1)-\left\lfloor \frac{(n-1+q-p)^2}{4}\right\rfloor-q\\
&=p-q+\left\lfloor \frac{(n+q-p)^2}{4}\right\rfloor-\left\lfloor \frac{(n-1+q-p)^2}{4}\right\rfloor
     \end{aligned}
\end{equation*}
We consider the parity of $n+q-p$.\\
If $n+q-p=2k$, then
$$\left\lfloor \frac{(n+q-p)^2}{4}\right\rfloor-\left\lfloor \frac{(n-1+q-p)^2}{4}\right\rfloor=k^2-\left\lfloor\frac{4k^2-4k+1}{4}\right\rfloor=k=\left\lfloor \frac{n+q-p}{2}\right\rfloor.$$
If $n+q-p=2k+1$, then
$$\left\lfloor \frac{(n+q-p)^2}{4}\right\rfloor-\left\lfloor \frac{(n-1+q-p)^2}{4}\right\rfloor=\left\lfloor\frac{4k^2+4k+1}{4}\right\rfloor-\frac{4k^2}{4}=k=\left\lfloor \frac{n+q-p}{2}\right\rfloor.$$
Hence,
$$d_{\overrightarrow{G}}(v)-q\geq p-q+\left\lfloor \frac{n+q-p}{2}\right\rfloor
=\left\lfloor \frac{n+p-q}{2}\right\rfloor.$$
Combining this with \eqref{a4}, we obtain
$$N(\overrightarrow{S}_{p,q},\overrightarrow{G})\geq \binom{\left\lfloor d_{\overrightarrow{G}}(v)-q\right\rfloor}{p}
\geq \binom{\left\lfloor \frac{n+p-q}{2}\right\rfloor}{p},$$
as desired.\\
CASE 2: $2\leq |Z|<4^p.$
Observe that
\begin{equation*}
     \begin{aligned}
F(n)-F(n-|Z|)&=\sum_{i=0}^{|Z|-1}(F(n-i)-F(n-i-1))\\
&=\sum_{i=0}^{|Z|-1}(p-1+\left\lfloor \frac{n+q-p-i}{2}\right\rfloor)\\
&=|Z|(p-1)+\sum_{i=0}^{|Z|-1}\left\lfloor \frac{n+q-p-i}{2}\right\rfloor\\
&\geq |Z|(p-1)+\sum_{i=0}^{|Z|-1}\left\lfloor \left(\frac{n+q-p}{2}-i\right)\right\rfloor\\
&=|Z|(p-1)+|Z|\left\lfloor\frac{n+q-p}{2}\right\rfloor-\frac{|Z|(|Z|-1)}{2}.
\end{aligned}
\end{equation*}
Combined with \eqref{a2}, 
\begin{equation*}
     \begin{aligned}
\frac{1}{|Z|}\sum_{v\in Z}\left(d_{\overrightarrow{G}}(v)-q\right)&\geq \frac{F(n)-F(n-|Z|)}{|Z|}+\frac{1}{|Z|}-q\\
&\geq p-1+\left\lfloor\frac{n+q-p}{2}\right\rfloor-\frac{|Z|-1}{2}+\frac{1}{|Z|}-q\\
&=\left\lfloor \frac{n+p-q}{2}\right\rfloor-\frac{|Z|+1}{2}+\frac{1}{|Z|}.
 \end{aligned}
\end{equation*}
Now apply \eqref{a4}:
\begin{equation*}
\begin{aligned}
N(\overrightarrow{S}_{p,q},\overrightarrow{G})
&\geq |Z|\binom{\left\lfloor\frac{1}{|Z|}\sum_{v\in Z}\big(d_{\overrightarrow{G}}(v)-q\big)\right\rfloor}{p}\\
&\geq |Z|\binom{\left\lfloor\left\lfloor \frac{n+p-q}{2}\right\rfloor-\frac{|Z|+1}{2}+\frac{1}{|Z|}\right\rfloor}{p}\\
&\geq |Z|\binom{\left\lfloor \frac{n+p-q}{2}\right\rfloor-4^p}{p}.
\end{aligned}
\end{equation*}
We estimate the binomial ratio: 
\begin{equation*}
     \begin{aligned}
\frac{\binom{\left\lfloor \frac{n+p-q}{2}\right\rfloor-4^p}{p}}{\binom{\left\lfloor\frac{n+p-q}{2}\right\rfloor}{p}}&=\prod_{i=0}^{p-1}\frac{\lfloor \frac{n+p-q}{2}\rfloor-4^p-i}{\lfloor \frac{n+p-q}{2}\rfloor-i}\\
&=\prod_{i=0}^{p-1}(1-\frac{4^p}{\lfloor \frac{n+p-q}{2}\rfloor-i})\\
&\geq \left(1-\frac{4^p}{\lfloor \frac{n+p-q}{2}\rfloor-p+1}\right)^p\\
&\geq \left(1-\frac{1}{2p}\right)^p\geq \frac{1}{2},
 \end{aligned}
\end{equation*}
where the last inequality follows from Bernoulli’s inequality.
Therefore, 
$$N\left(\overrightarrow{S}_{p,q},\overrightarrow{G}\right)\geq |Z|\binom{\lfloor \frac{n+p-q}{2}\rfloor-4^p}{p}\geq \frac{|Z|}{2}\binom{\lfloor \frac{n+p-q}{2}\rfloor}{p}\geq \binom{\lfloor \frac{n+p-q}{2}\rfloor}{p},$$
as desired.\\
CASE 3: $ 4^p\leq |Z|\leq n$. From \eqref{a2}, 
\begin{equation*}
     \begin{aligned}
\frac{1}{|Z|}\sum_{v\in Z}\left(d_{\overrightarrow{G}}(v)-q\right)&\geq \frac{F(n)-F(n-|Z|)}{|Z|}+\frac{1}{|Z|}-q\\
&=p-1+\frac{1}{|Z|}\left(\left\lfloor\frac{(n+q-p)^2}{4}\right\rfloor+1-\left\lfloor\frac{(n+q-p-|Z|)^2}{4}\right\rfloor\right)-q\\
&\geq p-1+\frac{1}{|Z|}\left(\frac{(n+q-p)^2}{4}-\frac{(n+q-p-|Z|)^2}{4}\right)-q\\
&\geq \frac{n+p-q}{2}-1-\frac{|Z|}{4}\\
&\geq \frac{n}{4}-\frac{q-p}{2}-1\geq \frac{n}{8}+p.
 \end{aligned}
\end{equation*}
By \eqref{a4}
\begin{equation*}
     \begin{aligned}
N\left(\overrightarrow{S}_{p,q},\overrightarrow{G}\right)&\geq |Z|\binom{\lfloor\frac{1}{|Z|}\sum_{v\in Z}(d_{\overrightarrow{G}}(v)-q)\rfloor}{p}\\
&\geq |Z|\binom{\lfloor\frac{n}{8}+p\rfloor}{p}\geq |Z|\frac{(\frac{n}{8})^p}{p!}\\
&\geq 4^p\frac{(\frac{n}{8})^p}{p!}=\frac{(\frac{n}{2})^p}{p!}\geq \binom{\lfloor\frac{n+p-q}{2}\rfloor}{p},
\end{aligned}
\end{equation*}
where the last inequality holds by $\lfloor\frac{n+p-q}{2}\rfloor\leq \frac{n}{2}$, as desired.\\
Therefore, for any oriented graph $\overrightarrow{G}$ satisfying the arc condition, we have
$$N(\overrightarrow{S}_{p,q},\overrightarrow{G})\geq \binom{\lfloor\frac{n+p-q}{2}\rfloor}{p}.$$
This finish the proof.
\end{proof}
\section{Proof of Theorem~\ref{S0p}}
{\bf Theorem~\ref{S0p}:} {\em    Fix a positive integers $q$, and let $\overrightarrow{G}$ be an $n$-vertex oriented graph, where $n$ is sufficiently large. If $|A(\overrightarrow{G})|=(q-1)n+k$ for some integer $k\geq 1$, then
$$N(\overrightarrow{S}_{0,q},\overrightarrow{G})\geq k.$$}
\begin{proof}
     We use $N(\overrightarrow{S}_{0,q},v)$ to count those  $\overrightarrow{S}_{0,q}$-copies in which vertex $v$ serves as the center. \\
    Clearly, $$N(\overrightarrow{S}_{0,q},v)=\binom{d^+(v)}{q}.$$\\
    Therefore, $$N(\overrightarrow{S}_{0,q},\overrightarrow{G})=\sum_{v\in V(\overrightarrow{G})}\binom{d^+(v)}{q}.$$We now prove the following claim:
    \begin{claim}\label{identity}
    For any integer $d\geq 0,$ 
        $$\binom{d}{q}\geq d-q+1.$$
    \end{claim}
    \begin{proof}[Proof of Claim \ref{identity}]
        When $d\leq q-1$, then $$\binom{d}{q}=0\geq d-q+1.$$
        We proceed by induction on $d$.
        When $d=q$, we have $$\binom{d}{q}=1=d-q+1.$$
        Assume the inequality holds for $d=k\geq q$, i.e. $\binom{k}{q}\geq k-q+1.$ \\
        By Pascal's identity,
    \begin{equation*}
     \begin{aligned}
 \binom{k+1}{q}&=\binom{k}{q}+\binom{k}{q-1}\\
 &\geq (k-q+1)+\binom{k}{q-1}\\
 &\geq (k-q+1)+1\\
 &=k-q+2.
  \end{aligned}
\end{equation*}
This finish the proof of Claim~\ref{identity}.
    \end{proof}

Therefore, 
\begin{equation*}
     \begin{aligned}
N(\overrightarrow{S}_{0,q},\overrightarrow{G})&=\sum_{v\in \overrightarrow{G}}\binom{d^+(v)}{q}\\
 &\geq\sum_{v\in{\overrightarrow{G}}}(d^+(v)-q+1)\\
 &=|A(\overrightarrow{G})|-n(q-1)\\
 &=(q-1)n+k-n(q-1)\\
 &=k,
  \end{aligned}
\end{equation*}
where the first inequality holds by Claim~\ref{identity}. This finish the proof.
\end{proof}
\vskip 1cm

\noindent {\bf Acknowledgement} \textnormal{This work was supported by National Natural Science Foundation of China under Grant Nos. 12371352 and 12261073 and Yongjiang Talent Introduction Programme of Ningbo under Grant No. 2021B-011-G.}

\vskip 3mm
\noindent {\bf Data availability}\textnormal{ No data was used for the research described in the article.}

\end{document}